\documentclass{article}
\usepackage{graphicx}

\newtheorem{thm}{Theorem}[section]
\newtheorem{lem}[thm]{Lemma}
\newtheorem{defn}[thm]{Definition}
\newtheorem{cor}[thm]{Corollary}
\newtheorem{prop}[thm]{Proposition}

\newtheorem{rem}[thm]{Remark}

\newfont{\bbc}{msbm10 scaled\magstep1}
\def\Bbb#1{\hbox{{\bbc #1}}}
\newfont{\bbbc}{msbm10 scaled\magstep0}

\newfont{\nam}{cmcsc10 scaled\magstep1}

\newcommand{\R}{{\bf R}}

\newcommand{\cJ}{{\cal J}}

\def\eproof{{\hfill{\vrule height5pt width3pt depth0pt}\medskip}}

\newcommand{\proof}{{\noindent\bf Proof.$\quad$}}

\newcommand{\Inv}{{\rm Inv\,}}

\newcommand{\Int}{{\rm int\,}}

\newcommand{\dir}{{\rm dir}}

\begin{document}
\begin{center}
{\Large Rigorous Numerics for Partial Differential
Equations: the Kuramoto-Sivashinsky equation}\\
\bigskip
by
\bigskip
\\
Piotr Zgliczy\'nski\footnote{Research supported in part by Polish KBN
grants 2P03A 021 15, 2 P03A 011 18
and NSF--NATO grant DGE--98--04459.}\\ 
 Jagiellonian University,
Institute of Mathematics, \\
Reymonta 4, 30-059 Krak\'ow, Poland \\ 

and \\
Center for Dynamical Systems and Nonlinear Studies\\
School of Mathematics\\
Georgia Institute of Technology\\
Atlanta, GA 30332 \\

e-mail: zgliczyn@im.uj.edu.pl, piotrz@math.gatech.edu

\bigskip

{\large and }

\bigskip

Konstantin Mischaikow\footnote{Research supported in part by
NSF  grant DMS-9805584.}\\
\ \\
Center for Dynamical Systems and Nonlinear Studies\\
School of Mathematics\\
Georgia Institute of Technology\\
Atlanta, GA 30332 \\

e-mail: mischaik@math.gatech.edu

\end{center}

\bigskip

\begin{abstract}
We present a new topological method for the study of the dynamics of 
dissipative PDE's. The method is based on the
concept of the self-consistent apriori bounds, which allows to
justify rigorously the Galerkin projection. As a result we obtain
a low-dimensional  system of ODE's subject to rigorously controlled
small perturbation from the neglected modes. To this ODE's
we apply the Conley index to obtain information about the dynamics of 
the PDE under consideration.

We applied the method to the Kuramoto-Sivashinsky equation 
\begin{displaymath}
    u_t=(u^2)_x - u_{xx} - \nu u_{xxxx}, \  u(x,t)=u(x +2\pi,t), \ u(x,t)=-u(-x,t) 
\end{displaymath}

We obtained a computer assisted proof the existence of the large number
fixed points for various  values of $\nu > 0$.
\end{abstract}

\bigskip


\vfill
\eject

\section{Introduction}

Even in the setting of infinite dimensional dynamics many of
the dynamical objects of interest are low dimensional, e.g.\
equilibria, periodic orbits, connecting orbits, horseshoes, etc.
In this paper we introduce techniques which, in principle, allow for
the rigorous verification of such solutions for a wide variety
of partial differential equations. Our approach is to combine rigorous
computer calculations with topological invariants to obtain
accurate existence statements. To demonstrate these techniques we
have chosen to study the the Kuramoto-Sivashinsky  equation \cite{KT,Si}
\begin{equation}
\label{eq:KS}
u_t = -\nu u_{xxxx} - u_{xx} + 2uu_x \quad\quad (t,x)\in
[0,\infty)\times (-\pi,\pi)
\end{equation}
 subject to
periodic and odd boundary conditions
\begin{equation}
\label{eq:ksbc}
u(t,-\pi) = u(t,\pi)\quad{\rm and}\quad
u(t,-x) = -u(t,x).
\end{equation}
The following theorem is a prototype for the results which can be
obtained.

\begin{thm}
\label{thm:proto}
Let $u(x) = \sum_{k=1}^{28} a_k\sin(kx)$ where the $a_k$ are given
in Table~\ref{table:coef}. Then, for $\nu=0.1$ there exists an
equilibrium
$u^*(x)$ for (\ref{eq:KS}) such that
\[
|| u^* - u||_{L^2} < 2.71547\times 10^{-13}\quad{\rm and}\quad
||u^*-u||_{C^0} < 2.06706\times 10^{-13}.
\]
\end{thm}

\begin{table}
\begin{tabular}{|c|c|c|}
\hline
$a_1 = 1.07934\times 10^{-37} $ & $a_2 =1.25665 $ & $a_3= -1.92524\times 10^{-37} $ \\
$a_4 = -0.559867 $ & $a_5=7.81863\times 10^{-38}  $ & $ a_6=0.0881138 $ \\
 $ a_7= -1.56596\times 10^{-38} $ & $ a_8=-0.0122945 $ & $a_9=2.54974\times 10^{-39} $ \\
$a_{10}= 0.00143504 $ & $a_{11}= -3.4963\times 10^{-40} $ & $ a_{12}=-0.000156065 $\\
$a_{13}= 4.35072\times 10^{-41} $ & $a_{14}= 1.59816\times 10^{-05} $ & $a_{15}=-5.02979\times 10^{-42} $ \\
$a_{16}= -1.57158\times 10^{-06} $&$a_{17}= 5.50953\times 10^{-43} $ & $a_{18}= 1.49677\times 10^{-07} $ \\
$a_{19}=-5.62586\times 10^{-44} $ & $a_{20}= -1.39049\times 10^{-08}$ & $a_{21}= -8.26547\times 10^{-45} $ \\
$a_{22}= 1.26591\times 10^{-09} $ & $ a_{23}=1.30584\times 10^{-43} $ & $a_{24}= -1.13347\times 10^{-10}$\\
$a_{25}=-9.46577\times 10^{-43} $ & $a_{26}= 1.0008\times 10^{-11} $ & $ a_{27}=1.1614\times 10^{-40} $ \\
$a_{28}= -8.73294\times 10^{-13}$ & &\\
   \hline
    \end{tabular}
\caption{Coefficients for the function $u(x)$.}
\label{table:coef}
\end{table}

\begin{figure}[h]
\begin{center}
\includegraphics[width=0.80\textwidth]{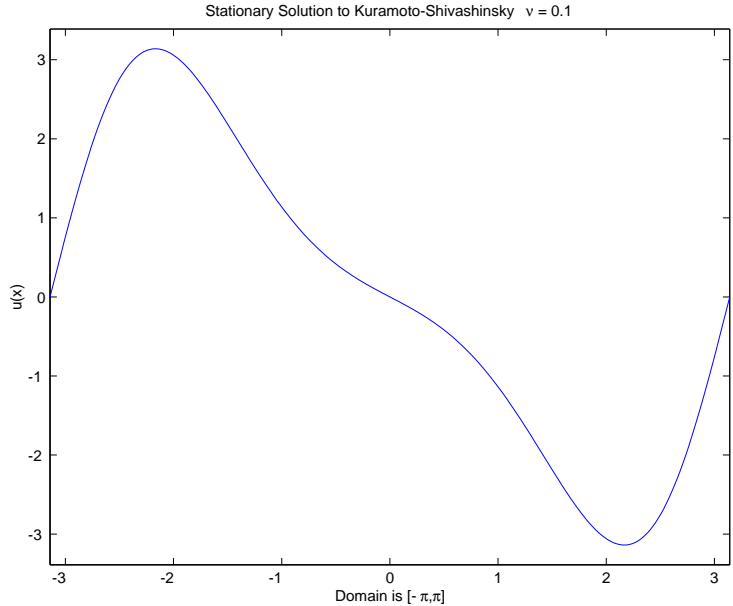}
\end{center}
\caption{The equilibrium solution $u^*(x)$ for $\nu=0.1$}
\label{fig:equil}
\end{figure}

Having stated this theorem we now try to put the result into the context of
the goals of our methods. To begin with it needs to be emphasized that the
computations which lead to this result are rigorous in the sense that
we have employed interval arithmetic to overcome all errors introduced by
the fact that we are using floating point arithmetic in our calculations.

As it will become clear in the later sections, this result is obtained by
studying the full partial differential equation rather than attempting to
solve a boundary value problem. While from the point of view of traditional
numerical analysis this approach may appear inefficient, it is an important
point. To be more precise, our method does not attempt to directly approximate
any particular solution to the partial differential equation. Rather we
essentially compute the Conley index of a compact region, called an
{\em isolating
neighborhood}, of
phase space.  The diameter of this region provides the error bounds stated in the
theorem. The index guarantees the existence of the equilibrium solution.

The Conley index theory is a far reaching topological generalization of
Morse theory. In particular,
this index can be used to prove the existence of periodic orbits,
connecting orbits, and chaotic dynamics
\cite{conley-cbms,smoller,salamon,cime,handbook}. It has been
numerically observed that for various parameter values (\ref{eq:KS})
contains these types of dynamical objects. In principle, combining
earlier rigorous numerical methods
\cite{mm-lorenz,mm-lorenz2,mms,az,gz,zg1,zg2} with
the techniques described in this paper and the above mentioned index
theorems will lead to rigorous proofs of the existence of periodic
orbits and even chaotic dynamics. However, we do not pursue these
more complicated structures in this paper for two reasons. First,
finding the appropriate isolating neighborhoods is more complicated
in these cases and our goal here is to emphasize the fundamental ideas
associated with the methods. The second, and more important point, is
that a straightforward application of the earlier numerical methods would
lead to large computations - which we believe can be avoided by
alternative methods (see for example \cite{wz}).
This latter point is currently being investigated.

Returning to our discussion of Theorem~\ref{thm:proto}, an obvious
question concerns the stability of $u^*$. For this we have no definitive
answer. As was indicated before our method does not directly approximate
$u^*$ and therefore we do not obtain uniqueness results or hyperbolicity
results. On the other hand, being a generalization of the Morse index
the Conley index does contain some information about the stability
or instability of the dynamics in the isolating neighborhood.
Thus, what
can be asserted is the following. Assume that
$u^*$ is a hyperbolic fixed point, i.e.\ all eigenvalues have
nonzero real part, and that $u^*$ is the only solution which remains within
either the $L^2$ or $C^0$ bounds of $u$ for all time, then
$u^*$ has exactly two unstable eigenvalues, i.e.\ its unstable manifold is
two dimensional. We hope to treat this problem in a subsequent paper.

It should be mentioned that even though we are doing the computations via
an approximation of the full partial differential equation, we never
integrate the equations. Rather, as will be made clear in Section 2
the computations are reduced to solving a set of inequalities. It is
for this reason that we are able to get such sharp bounds on the
equilibria. As the following theorem demonstrates we can, in fact
obtain bounds on the level of the floating point accuracy.

\begin{thm}
\label{thm:stable}
One can compute a sequence $a_1,a_2,\dots,a_{30}$ and
the function  $u(x) = \sum_{k=1}^{30} a_k\sin(kx)$, such that
 for $\nu=0.75$ there exists an equilibrium
$u^*(x)$ for (\ref{eq:KS}) such that
$
|| u^* - u||_{L^2} < 1.26281\times 10^{-15}$ and
$||u^*-u||_{C^0} < 9.57396\times 10^{-16}$.
\end{thm}

By now it is a well demonstrated principle that the asymptotic
behavior of a wide variety of infinite dimensional dynamical systems
is finite dimensional \cite{hmo,hale,temam}. The Kuramoto-Sivashinsky
equation (\ref{eq:KS}) is a particularly well studied example of
such a system \cite{FNST}. In fact, it is known that (\ref{eq:KS}) posses an
inertial manifold and therefore, that there exists a family of
ordinary differential equations that exactly describes the asymptotic
dynamics. Unfortunately, the estimates for the dimension of
these manifolds make them impractical for our purposes \cite{JRT}.

We mention these methods to emphasize that our approach does not
directly make use of any of these results. What appears to be essential for
our techniques is that the spectrum of the linear operator for the evolution
equation is not clustered near the imaginary axis. This is in contrast to
the inertial manifold techniques which strongly rely on gap conditions
of the spectrum or cone conditions from the flow. Our approach is
to use the computer to restrict our attention to that portion of phase
space in which the desired dynamics (for this paper the fixed points)
occur. Obviously, by restricting the phase space one can get much
better estimates. This sets up a loop by which one can continuously
improve the estimates until the desired bounds are reached.

Our analysis of the fixed points for (\ref{eq:KS}) was motivated in
part by the work of Jolly, Keverkidis, and Titi \cite{jkt}. In
particular, using a 12 mode traditional Galerkin approximation, they
produced a bifurcation diagram for $\nu \in (0.057,\infty)$. We used their
reported solutions to test our methods. In particular, as is indicated
below we were able to find and prove the existence of an equilibrium
point on each of their stable branches. Unfortunately, we used a
fairly primitive search procedure and therefore missed a few unstable
branches. Our expectation is that by combining our methods with a
continuation package, one could produce a rigorous bifurcation
diagram with fairly precise bounds in a computationally inexpensive
manner.

Below we include some of the  steady states we found.

\begin{itemize}
\item{$\nu=0.5$}.  Two stable unimodal fixed points

\item{$\nu=0.3$}.
      Two stable unimodal fixed points

\item{$\nu=0.127$, $\nu=0.125$}.  A stable and unstable
    bimodal fixed point.
    Negative branch is stable, positive one is unstable
    with apparently two-dimensional unstable manifold.

    Our primitive search procedure did not find a solution on
the bi-tri branch.

\item{$\nu=0.1$}
   An bimodal stable and unstable (2 unstable directions)
   and two unstable trimodal fixed points (both with 1-dimensional unstable
   manifold)

  We did not find an unstable branch connecting bi-tri branch with quadrimodal
  branch.

\item{$\nu=0.08$}
    A bimodal stable (neg. branch) and unstable (2 unstable directions) fixed
    points.
    A pair of stable fixed points close to $R_3 t_2$ (see \cite{jkt}).
    A pair of unstable trimodal fixed points (1 unstable direction).

We did not find  a branch connecting bi-tri and quadimodal branches.

\item{$\nu=0.0666..$, $\nu=0.063$}
   Two unstable bimodal points, two stable trimodal points
   and two stable solutions apparently belonging to the {\em giant }
   branch.

  We are lacking two unstable branches which are present in \cite{jkt}.

\item{$\nu=0.062$},
   Two stable trimodal points and two stable points from giant branch

\item{$\nu=0.045$},
  Two stable points from giant branch and pairs of unstable tri- and
  quadrimodal fixed points

\item{$\nu=0.04$},
Two stable giant fixed points.
Two stable quadrimodal fixed points.
Two unstable trimodal points.

\item{$\nu=0.029$}
  Two unstable quadrimodal points.

\end{itemize}

\section{The Method}

Our method begins with the reduction of the full dynamical system to
a lower dimensional system which can be studied numerical. In
particular, we begin with a nonlinear evolution equation in a
Hilbert space $H$ ($L^2$ in our treatment of Kuramoto-Sivashinsky)
of the form
\begin{equation}
\label{eq:pde}
\frac{du}{dt} = F(u)
\end{equation}
where domain of $F$ is  dense in $H$. Furthermore, we assume that
$\{ \phi_i\mid i=0,1,\ldots\}$ forms a complete orthogonal
basis for $H$.

In the case of the Kuramoto-Sivashinsky equation $F(u)=Lu + B(u,u)$,
where $L$ is a linear part and $B$ is a nonlinear part, the functions
$\{\varphi_i\}$ are chosen to be eigenvalues of $L$.

Fix $m \in {\Bbb N}$ and  let
\[
P=P_m:H\to X_m =X
\]
be the orthogonal projection
from $H$ onto the finite
dimensional subspace span\-ned by $\{ \phi_1,\phi_2,\dots,\phi_m\}$.
Let
\[
Q=Q_m:= I-P:H\to Y=Y_m
\]
be the complementary orthogonal projection.
Finally, let
\[
A_k : H\to {\Bbb R}
\]
 be the orthogonal projection onto
the subspace generated by $\phi_k$.

Given
$u\in H$, let $Pu = p$ and $Qu =q$. Equation (\ref{eq:pde}) can be
rewritten as
\begin{eqnarray}
\label{eq:galode}
\frac{dp}{dt} & = & PF(p,q) \\
\label{eq:comode}
\frac{dq}{dt} & = & QF(p,q)
\end{eqnarray}

The strategy adopted is fairly simple: study the
dynamics of
the low dimensional Galerkin projection (\ref{eq:galode}) to draw
conclusions about the dynamics of (\ref{eq:pde}). Before turning to the
precise conditions, consider the following  heuristic description of our
approach.

Let
$W\subset X=X_m$. For
$j > m$,  let $W_j\subset X_j$ such that $P(P_j^{-1}(W_j))=W$,
(i.e. $W_j=W \oplus (I-P) W_j$).
Similarly, let $V\subset Y$ and set $V_j= Q_j(V)$. Furthermore, given
$q_j\in V_j$ assume that
$\lim_{j\to\infty} ||q_j|| = 0$.
Our only knowledge concerning the higher order modes or ``tails'' of
the solutions is that they project into $V$. This implies that our
knowledge of the vector field is reduced to the following
differential inclusion
\[
\frac{dp}{dt}  \in   PF(p,V)
\]
where $p\in W$.
Numerical calculations on this equation are used to find topological
invariants (the Conley index, the fixed point index) which guarantee
the existence of specific dynamics, e.g. fixed points, periodic orbits,
symbolic dynamics, positive entropy, etc. It is simultaneously
argued that the topological invariant is the same for any
Galerkin system  of the form
\[
\frac{dp_j}{dt}  \in   PF(p_j,V_j)
\]
where $p_j\in W_j$. Thus, the same dynamical object exists
for each sufficiently high Galerkin approximation. Finally,
it is shown that the limit of these objects leads to the desired
dynamics for the full system (\ref{eq:pde}).

\subsection{Self-consistent Bounds}

As one might expect the orthomormal basis $\{\phi_i \}$ and the sets
$W$ and $V$ must be chosen with
care. The first issue that needs to be dealt with is analytic
in nature - solutions to the  ordinary differential
equations must approximate solutions of the partial differential
equation. This
leads to the following definition.

\begin{defn}
\label{defn:selfconsistant}
{\em
Let $m,M\in{\Bbb N}$ with $m\leq M$. A compact set $W\subset X_m$
and a sequence
of pairs
$\{ a_k^\pm\in{\Bbb R}\mid a_k^-<a_k^+,\ k\in{\Bbb N}\}$ form {\em
self-consistent apriori bounds} for (\ref{eq:pde}) if the following
conditions are satisfied:
\begin{description}
\item[C1] For $k>M$, $a_k^- < 0< a_k^+$.
\item[C2] Let $\hat{a}_k := \max |a_k^\pm |$ and set
$\hat{u} = \sum_{k=0}^\infty \hat{a}_k\phi_k$. Then,
$\hat{u}\in H$.
In particular, $||\hat{u}||<\infty$.
\item[C3] The function $u\mapsto F(u)$ is continuous on
\[
W\oplus\prod_{k=m+1}^\infty [a_k^-,a_k^+]\subset H.
\]
\end{description}
}
\end{defn}

In practice $W\subset \prod_{k=1}^m [a_k^-,a_k^+]$. Given
self-consistent apriori bounds $W$
and $\{ a_k^\pm\}$, let
\[
V:=\prod_{k=m+1}^\infty [a_k^-,a_k^+]\subset
Y_m.
\]
Our goal is to numerically solve (\ref{eq:galode}) on $W$
and draw conclusions about the dynamics of (\ref{eq:pde})
on the set $W\oplus V\subset H$. To do this we will make use of
the following results, the first two of them are obvious.

\begin{lem}
\label{lem:compact}
Given
self-consistent apriori bounds $W$
and $\{ a_k^\pm\}$, $W\oplus V$ is a compact subset of $H$.
\end{lem}

\begin{lem}
\label{lem:compact2}
Given self-consistent apriori bounds $W$
and $\{ a_k^\pm\}$, $W\oplus V$, then
\begin{displaymath}
  \lim_{n\to \infty} Q_n(F(u)) =0, \quad \mbox{uniformly for
   $u \in W \oplus V$}
\end{displaymath}
\end{lem}

\begin{prop}
\label{prop:ladder}
Let $W$
and $\{ a_k^\pm\}$ be self-consistent bounds for (\ref{eq:pde}).
A function $a:[0,T]\to W\oplus V$ given by
\[
a(t) := \sum_{k=0}^\infty a_k(t)\phi_k
\]
is a solution to (\ref{eq:pde}), if and only if, for each
$k\in{\Bbb N}$ and all $t\in[0,T]$
\begin{equation}
\label{eq:ladder}
\frac{da_k}{dt} =A_kF(a).
\end{equation}
\end{prop}

\proof ($\Rightarrow$) This direction follows directly from the
projection of (\ref{eq:pde}) onto each of the basis elements.

($\Leftarrow$) Assume that (\ref{eq:ladder}) is satisfied for
each $k\in{\Bbb N}$ and all $t\in [0,T]$. Let
\[
a(t) := \sum_{k=0}^\infty a_k(t)\phi_k \in H
\]
First observe that from {\bf C3} it follows immediately that
$\sum_{k=1}^{\infty}\frac{d a_k}{dt} \phi_k = F(a) \in H$.

It needs to be shown that
\[
\frac{da}{dt} =\lim_{h\to 0}\frac{a(t+h) -a(t)}{h} = \sum_{k=0}^\infty
\frac{da_k}{dt} \phi_k.
\]
This is equivalent to showing that
\[
\lim_{h\to 0}  \left|\frac{1}{h}\sum_{k=1}^\infty (a_k(t+h) -a_k(t))\phi_k -
\sum_{k=1}^{\infty}
\frac{d a_k}{dt} \phi_k \right| =0
\]
for all $t\in [0,T]$.

Fix $h>0$,  then for any $n\in{\Bbb N}$
\begin{eqnarray*}
\left|\frac{1}{h}\sum_{k=1}^\infty (a_k(t+h) -a_k(t))\phi_k -
\sum_{k=1}^{\infty}
\frac{d a_k}{dt} \phi_k \right| \leq \\
\left| \frac{1}{h}\sum_{k=1}^n (a_k(t+h) -a_k(t)) \phi_k - \sum_{k=1}^{n}
\frac{d a_k}{dt} \phi_k \right| \\
  + \left|\frac{1}{h}\sum_{k=n+1}^{\infty} (a_k(t+h)
-a_k(t)) \phi_k\right|+
  \left| \sum_{k=n+1}^{\infty} \frac{d a_k}{dt} \phi_k \right|
\end{eqnarray*}

We will estimate the three terms on the right hand side separately.
From lemma \ref{lem:compact2} it follows for a given $\epsilon
>0$ there exists $n_0$ such that $n>n_0$ implies
\[ \left| \sum_{k=n+1}^{\infty} \frac{d a_k}{dt} \phi_k \right| =
\left| Q_n(F(a)) \right|  < \epsilon/3.
\]
From now on fix $n>n_0$. Again lemma \ref{lem:compact2} and the mean value
theorem implies
\begin{eqnarray*}
\left|\sum_{k=n+1}^{\infty} \frac{1}{h}(a_k(t+h) -a_k(t))\phi_k\right| &=&
\left|\sum_{k=n+1}^{\infty} \frac{d a_k}{d t}(t+\theta_k h) \phi_k\right| \\
&=&
   \left| Q_n(F(a(t+\theta_k h)) \right|   < \epsilon/3 .
\end{eqnarray*}
Finally, for $h$ sufficiently small,
\[
\left|\frac{1}{h}\sum_{k=1}^n (a_k(t+h) -a_k(t))\phi_k - \sum_{k=1}^{n}
\frac{d a_k}{dt}\phi_k \right|
  < \epsilon/3
\]
and hence the desired limit is obtained.
\eproof

\subsection{Conley Index}

Proposition~\ref{prop:ladder} indicates that given self-consistent
apriori bounds $W$ and $\{ a^\pm_k\}$, finite time solutions to
(\ref{eq:ladder}) are solutions to the full partial differential
equation. Thus, the goal of this paper is to find solutions to
(\ref{eq:ladder}). Of course, numerically one can only study
(\ref{eq:galode}) restricted to
$W$ and then argue that the resulting numerical solution is an
approximation to a solution to (\ref{eq:pde}). Hence, rather than
attempting to approximate specific  trajectories in $W$ directly, the
objective is to
compute a Conley index for  (\ref{eq:galode}) and then show that this
index information is sufficient to guarantee a solution for
(\ref{eq:pde}).

In order to describe this index the following definitions are needed.
Let $\varphi :{\Bbb R}\times {\Bbb R}^m\to {\Bbb R}^m$ be a
continuous flow generated by a differential equation $\dot{z}=f(z)$.

\begin{defn}
\label{defn:isoneighbor}
{\em
A compact set  $N\subset {\Bbb R}^n$ is an {\em
isolating neighborhood} if
\[
\Inv(N,\varphi) := \{ z\in N\mid \varphi({\Bbb R},z)\subset N\}\subset
\Int N.
\]
If, in addition, for any $z\in\partial N$, there exists $t_z>0$
such that
\begin{equation}
\label{eq:exitentrance}
\varphi((0,t_z),z)\cap N=\emptyset\quad{\rm or}\quad
\varphi((-t_z,0),z)\cap N=\emptyset,
\end{equation}
then $N$ is an {\em isolating
block}. Given an isolating neighborhood
$N$, the
associated maximal invariant set $\Inv(N,\varphi)$ is an
{\em isolated invariant set}.
}
\end{defn}

The easiest way to verify the existence of an isolating block
is through
local sections.

\begin{defn}
\label{defn:section}
{\em
$\Xi \subset {\Bbb R}^n$ is a
{\em local section} for $\varphi$ if for
some $\epsilon >0$
\begin{equation}
\varphi:(-\epsilon,\epsilon)\times\Xi \to
\varphi((-\epsilon,\epsilon),\Xi)
\end{equation}
is a homeomorphism and $\varphi((-\epsilon,\epsilon),\Xi)$ is an
open subset of ${\Bbb R}^n$.
}
\end{defn}

A special form of local section is  a hypersurface which is
transverse to the flow. More formally,
let $\Xi\subset {\Bbb R}^n$ be an $n-1$ dimensional manifold with
normal vector $\mu(z)$ at $z\in\Xi$. $\Xi$ is a  local section if for
each
$z\in \Xi$,
\begin{equation}
\label{eq:section}
\mu(z)\cdot f(z)\neq 0.
\end{equation}

It is straightforward to check that $N$ is an isolating block
if $\partial N$ can be written as the union of the closure of
local
sections with the property that (\ref{eq:exitentrance}) is satisfied at
every point in the intersection of the closure  of the sections.

In this paper the focus is both on proving the existence of equilibria
and providing tight bounds on the location of the equilibria. To do
this requires have good isolating blocks. With this in mind consider
the linear ordinary differential equation
\begin{equation}
\label{eq:linearode}
\dot{z}= Bz,\quad z\in{\Bbb R}^n.
\end{equation}
Assume that the origin is a hyperbolic fixed point. Without loss of
generality it can be assumed that
$B$ is in Jordan normal form. Generically, to each real eigenvalue there
is associated a 1-dimensional eigenspace and to each pair of complex
conjugate eigenvalues there is an associated 2-dimensional eigenspace.
Thus,   ${\Bbb R}^n$ can be decomposed into the product of
eigenspaces, i.e.\
\[
{\Bbb R}^n= V_1\times V_2\times\cdots\times V_k
\]
where $V_i$ is either
${\Bbb R}$ or ${\Bbb R}^2$. In what follows we will use the following
notation, $z_i\in V_i$, $i=1,\ldots,k$, and if $V_i\cong {\Bbb R}^2$,
then $z_i=(x_i,y_i)$.

Our interest is not on the dynamics of (\ref{eq:linearode}) on the
entire phase space, but rather on a prescribed compact subset. Since
our goal is to understand the equilibria of (\ref{eq:linearode}) consider
a neighborhood of the origin,
\[
N=I_1\times I_2\times \cdots \times I_k
\]
where
\[
I_i := \cases{ [b_i^-,b_i^+],\ b_i^-<0 <b_i^+ & if $V_i \cong{\Bbb R}$,
\cr
\{ (x_i,y_i)\in V_i\mid \sqrt{x_i^2+y_i^2}\leq b_i,\ b_i>0\} &
if $V_i \cong{\Bbb R}^2$. \cr}
\]

The following result is obvious, but to  make a point crucial to the
results of this paper we will provide the proof.

\begin{lem}
The compact set $N$
is an isolating block for (\ref{eq:linearode}).
\end{lem}

\proof Since $B$ is in Jordan normal form the system decouples
according to the decomposition ${\Bbb R}^n= V_1\times
V_2\times
\cdots\times V_k$.

If $V_i={\Bbb R}$, then (\ref{eq:linearode})
reduces to $\dot{z}_i=\lambda_i z_i$. Since $B$ is hyperbolic
$\lambda_i\neq 0$, and hence at
$z_i = b^\pm_i$ (\ref{eq:section}) becomes
\[
\lambda_i b^\pm_i \neq 0.
\]

If $V_i={\Bbb R}^2$, then (\ref{eq:linearode})
reduces to
\begin{eqnarray*}
\dot{x}_i & = & \alpha_i x_i +\beta_i y_i \\
\dot{y}_i & = & -\beta_i x_i + \alpha_i y_i
\end{eqnarray*}
where by hyperbolicity $\alpha_i \neq 0$. So again
for $\sqrt{x_i^2+y_i^2}= b_i$ (\ref{eq:section})
becomes
\[
(x_i,y_i) \cdot (\alpha_i x_i +\beta_i y_i,-\beta_i x_i + \alpha_i y_i)^t
= \alpha_i b_i^2 \neq 0.
\]
\eproof

To see why this trivial argument is of importance, consider the more
interesting example of
\begin{equation}
\label{eq:nonlinearode}
\dot{z}=Bz +f(z)+E(z)
\end{equation}
where
$f:{\Bbb R}^n\to{\Bbb R}^n$ is $o(||z||^2)$ at $0$ and $E$ represents
a known bounded error. In our situation $E$ arises from numerical
errors and approximations. More
precisely, we assume that there  are  known constants
$c_i$ such that
\[
\sup_{z\in N}||E_i(z)||\leq c_i.
\]

Observe that a sufficient condition for $N$ to be an isolating
block for (\ref{eq:nonlinearode}) is the following: for each $i$ such that
$V_i={\Bbb R}$,
\begin{equation}
\label{eq:isolatinginequality1}
\lambda_i b_i^\pm + f_i(z) + E_i(z)
\end{equation}
has the same sign as $\lambda_i b^\pm_i$ over the set
$\{ z\in N\mid z_i=b_i^\pm\}$, and for each $i$ such that
$V_i={\Bbb R}^2$,
\begin{equation}
\label{eq:isolatinginequality2}
(x_i,y_i) \cdot (\alpha x_i +\beta y_i +f_{i_1}(x)+E_{i_1}(x),-\beta x_i +
\alpha y_i+f_{i_2}(x)+E_{i_2}(x))^t
\end{equation}
has the same sign as $\alpha_i$ over the set $\{ z\in N
\mid \sqrt{x_i^2+y_i^2}= b_i\}$.

For the
linear case the eigenvalues of $B$ are assumed to be known and
the $b_i^\pm$ can be chosen arbitrarily. Therefore, one can also interpret
(\ref{eq:isolatinginequality1}) and (\ref{eq:isolatinginequality2}) as
providing a set of inequalities that if simultaneously solved
for $b_i^\pm$ provide
an isolating block even in the context of numerical errors and
approximations. In particular, finding isolating blocks need not involve
numerically solving the ordinary differential equation.

In itself the knowledge that $N$ is an isolating block does not imply
anything about $\Inv(N,\varphi)$. To gain information concerning
the isolated invariant set we will make use of the Conley index.
For our purposes we need only a very small portion of the index
theory and so we give a minimal operational definition (see
\cite{conley-cbms,smoller,salamon,cime,handbook} for further
information).

\begin{defn}
{\em
Let $N$ be an isolating block and let $\partial N = L^+\cup L^-$
where $L^\pm$ are closed sets. Furthermore, assume that
$z\in L^-$ implies that
\[
\varphi((0,\epsilon),z)\cap N =\emptyset
\]
for a sufficiently small $\epsilon >0$. Similarly, assume that
if $z\in \cap L^+$, then
\[
\varphi((-\epsilon,0),z)\cap N =\emptyset
\]
for a sufficiently small $\epsilon >0$.
 The {\em Conley index}
of $S=\Inv(N,\varphi)$ is
\[
CH_*(S) := H_*(N,L).
\]
}
\end{defn}

No knowledge of relative
homology groups is required for the applications described in this
paper. The following theorem gives a formula for the index of
a hyperbolic fixed point.

\begin{prop}
\label{prop:indexcomp}
Let $q$ be the number of eigenvalues of $B$ with positive real
part. Assume that for all $i=1,\ldots ,k$ either the condition
associated with
(\ref{eq:isolatinginequality1}) or the
condition associated with (\ref{eq:isolatinginequality2}) are satisfied .
Then,
\[
CH_j(\Inv(N,\varphi))\cong\cases{ {\Bbb Z} & if $j=q$\cr
0 & otherwise.}
\]
\end{prop}

The following theorem due to McCord \cite[Corollary 5.9]{mccord},
indicates that if  the Conley index is a that of
Proposition~\ref{prop:indexcomp}, then there exists a  fixed point
in $N$.

\begin{thm}
\label{thm:fixedptindex}
If the Conley index has the form
\[
CH_j(\Inv(N,\varphi)) \cong\cases{ {\Bbb Z} & if $j=q$\cr
0 & otherwise,}
\]
for some $q$, then $N$ contains a fixed point.
\end{thm}

To indicate how these index ideas will be used in this paper let us
return to the system (\ref{eq:nonlinearode}). Observe that the only
assumption on the error term $E$ was that it is bounded, therefore,
it is no longer apriori true that the origin is a fixed point
or that there even exists a fixed point to (\ref{eq:nonlinearode}).
On the other hand the sets $N$ and $L$ remain unchanged. Therefore,
the Conley index implies the existence of the fixed point.

\subsection{A Singular Perturbation Result}

As was indicated in the previous section, it is possible to
find an isolating block for a finite dimensional ordinary
differential equation about a fixed
point by solving an appropriate set of inequalities. However,
to do this requires a good estimate of the location of the
fixed point, knowledge of the eigenvalues, the ability to evaluate the
nonlinear terms, and estimates of associated errors. Therefore, the
dimension to which one can hope to apply this procedure is obviously
limited. In this section we will describe a singular perturbation result
which allows one to ``lift'' the index computations of the previous
sections to arbitrary dimensions.

The definition of self-consistent bounds related individual solutions
of the
infinite family of ordinary differential equations to solutions in
partial differential equation. We now need to extend this definition
in order to know that the index computations we perform for the
finite dimensional approximation  have implications for the
partial differential equation.

\begin{defn}
{\em
Let $m,M\in{\Bbb N}$ with $m\leq M$. A pair of compact sets
$
N\subset W\subset X_m
$
and a sequence
of pairs
$\{ a_k^\pm\in{\Bbb R}\mid a_k^-<a_k^+,\ k\in{\Bbb N}\}$ are {\em
topologically self-consistent}  if $W$ and $\{ a_k^\pm\}$ are
self-consistent apriori bounds and the following conditions are satisfied.
\begin{description}
\item[C4] Let $u\in W\oplus\prod_{k=m+1}^\infty [a_k^-,a_k^+]$.
Then, for $k>m$
\begin{equation}
\label{eq:general}
A_k u = a_k^\pm \quad\Rightarrow\quad  A_kF(u) \neq 0.
\end{equation}
\item[C5] $N$ is an isolating block for (\ref{eq:galode})
for all $q\in \prod_{k>m}[a_k^-,a_k^+]$.
\end{description}
}
\end{defn}

For Kuramoto-Sivashinsky we will make use of the following stricter
form of {\bf C4}.

\begin{description}
\item[C4a] Let $u\in W\oplus\prod_{k=m+1}^\infty [a_k^-,a_k^+]$.
Then, for $k>m$
\begin{eqnarray}
\label{eq:up}
A_k u = a_k^+ &\Rightarrow & A_kF(u)) < 0 \\
\label{eq:down}
A_k u = a_k^- &\Rightarrow & A_kF(u)) >0.
\end{eqnarray}
\end{description}

Using the line of reasoning that was described in the analysis of
(\ref{eq:nonlinearode}), condition {\bf C5} can be replaced by
the following assumption.

\begin{description}
\item[C5a]  Let $N$ be an isolating block for (\ref{eq:galode}).
Let
$\nu^\pm(p)$ be the outward normal at
$p\in  L^\pm$. If $u\in W\oplus\prod_{k=m+1}^\infty
[a_k^-,a_k^+]$ such that $Pu\in L^\pm$, then
\[
PF(u)\cdot\nu^+(Pu) <0 \quad\quad PF(u)\cdot\nu^-(Pu) >0.
\]
\end{description}

We shall now discuss two singular perturbation results. The first
allows one to lift isolating blocks.
\begin{thm}
\label{thm:isoblock}
Let $m,M\in{\Bbb N}$ with $m\leq M$. Assume
$
N\subset W\subset X_m
$
and the sequence
of pairs
$\{ a_k^\pm\in{\Bbb R}\mid a_k^-<a_k^+,\ k\in{\Bbb N}\}$ are
topologically self-consistent. Fix an integer $r >m$. Then
for any $q=\sum_{k=r+1}^\infty q_k \phi_k $, such that
$q \in \Pi_{k=r+1}^\infty [a_k^-,a_k^+]$ and $q_k =0$ for $k>M$ the set
\[
\tilde{N}:= N\times [a_{m+1}^-,a_{m+1}^+]\times [a_{m+2}^-,a_{m+2}^+]
\times \cdots\times [a_{r}^-,a_{r}^+]
\]
is an isolating block for the system of equations
\begin{equation}
\dot{x}_k = A_kF(\sum_{i=1}^k x_i \phi_i + q)\quad\quad k=1,\ldots ,r
\label{eq:isoblock}
\end{equation}
where $x\in\R^r$.
\end{thm}

\proof Let $u=(w,v)\in W\oplus \prod_{k=m+1}^r[a_k^-,a_k^+]$. From
{\bf C1} it follows that $u + q \in  W \oplus \Pi_{k=m+1}^\infty[a_k^-,a_k^+]$.
If $u\in\partial \tilde{N}$, then either $w$ is in $\partial N$ or
$v$ is in $\partial \prod_{k=m+1}^r[a_k^-,a_k^+]$. In the first case
{\bf C5} forces the vector field to be transverse at the
boundary. In the second case transversality follows from
{\bf C4}.
\eproof

%\begin{rem}
%The proof of the above theorem is the only place where we use
%assumption {\bf C1}.
%\end{rem}

\begin{rem}
For $r>M$  equations (\ref{eq:isoblock}) are the Galerkin projection
of $\dot{u}=F(u)$.
\end{rem}

The direction of the vector field influences the index computation.
With this in mind define
\[
\dir(k) := \cases{
-1 & if $A_k u = a_k^+\ \Rightarrow\ A_kF(u)) < 0$ and\cr
&\ \ \ $A_k u = a_k^-\ \Rightarrow\  A_kF(u)) >0 $ \cr
0 & if $A_k u = a_k^+\ \Rightarrow\  A_kF(u)) < 0$ and\cr
&\ \ \ $A_k u = a_k^-\ \Rightarrow\  A_kF(u)) <0 $ \cr
0 & if $A_k u = a_k^+\ \Rightarrow\  A_kF(u)) > 0$ and\cr
&\ \ \ $A_k u = a_k^-\ \Rightarrow\  A_kF(u)) >0 $ \cr
1 & if $A_k u = a_k^+\ \Rightarrow\  A_kF(u)) > 0$ and\cr
&\ \ \ $A_k u = a_k^-\ \Rightarrow\  A_kF(u)) <0 $ \cr}
\]

\begin{thm}
\label{thm:indexcomp}
Let $m,M\in{\Bbb N}$ with $m\leq M$. Assume
$
N\subset W\subset X_m
$
and the sequence
of pairs
$\{ a_k^\pm\in{\Bbb R}\mid a_k^-<a_k^+,\ k\in{\Bbb N}\}$ are
topologically self-consistent. Fix an integer $r >m$.  Let
$q=\sum_{k=r+1}^\infty q_k \phi_k $, such that
$q \in \Pi_{k=r+1}^\infty [a_k^-,a_k^+]$ and $q_k =0$ for $k>M$. Let
\[
\tilde{N}:= N\times [a_{m+1}^-,a_{m+1}^+]\times [a_{m+2}^-,a_{m+2}^+]
\times \cdots\times [a_{r}^-,a_{r}^+]
\]
Consider the dynamical system induced by (\ref{eq:isoblock}).

If for some $j\in\{ m+1,\ldots ,r\}$, $\dir(j)=0$, then
\[
CH_*(\Inv(\tilde N)) =0.
\]
Assume that for all $j\in\{ m+1,\ldots ,r\}$, $\dir(j)\neq 0$,
and let $d$ be the number of $j\in\{ m+1,\ldots ,r\}$ such that
$\dir(j) = 1$, then
\begin{equation}
\label{eq:indexshift}
CH_{s+d}(\Inv (\tilde{N})) \cong CH_s(\Inv (N)).
\end{equation}
\end{thm}

\proof By Theorem~\ref{thm:isoblock}, $\tilde{N}$ is an isolating
block.

We will present the proof of the second part of the theorem, only.

Assume that for all $j\in\{ m+1,\ldots ,r\}$, $\dir(j)\neq 0$.
 Let $\cJ:=\{ j\mid m<j\leq r,\ \dir(j)=1\}$. Set
\begin{eqnarray*}
& \tilde{L}^- := & \left(L^-\times \prod_{k=m+1}^r [a_k^-,a_k^+]\right)\cup \\
&  & \bigcup_{j\in\cJ}\left(N\times \prod_{k=m+1}^{j-1} [a_k^-,a_k^+]\times
\{ a_j^\pm\}\times \prod_{k=j+1}^r [a_k^-,a_k^+]\right)
\end{eqnarray*}
and
\begin{eqnarray*}
& \tilde{L}^+ := & \left(L^+\times \prod_{k=m+1}^r [a_k^-,a_k^+]\right)\cup \\
& & \bigcup_{j\not\in\cJ}\left(N\times \prod_{k=m+1}^{j-1} [a_k^-,a_k^+]\times
\{ a_j^\pm\}\times \prod_{k=j+1}^r [a_k^-,a_k^+]\right)
\end{eqnarray*}
Let $\varphi:{\Bbb R}\times {\Bbb R}^r\to {\Bbb R}^r$, be any
flow generated by
\[
\dot{a}_k = A_kF(u + q)\quad k=1,\ldots,r
\]
where $a_k = A_k u$ and $u\in W\oplus \prod_{k=m+1}^r[a_k^-,a_k^+]$.

Clearly, if
$P_ru\in L^-$, then
$\varphi((0,\epsilon),P_ru)\notin \tilde{N}$ for small $\epsilon >0$.
Similarly,  if $P_ru\in L^+$, then
$\varphi((-\epsilon,0),P_ru)\notin \tilde{N}$ for small $\epsilon >0$.

Let $u=(w,v)\in W\oplus \prod_{k=m+1}^r[a_k^-,a_k^+]$. If
$u\in\partial \tilde{N}$, then either $w$ is in $\partial N$ or
$v$ is in $\partial \prod_{k=m+1}^r[a_k^-,a_k^+]$. Therefore,
$\partial N = \tilde{L}^+\cup \tilde{L}^-$.

Thus, the Conley index of $\Inv(\tilde{N})$ is given by
\[
CH_*(\Inv(\tilde{N})\cong H_*(\tilde{N},\tilde{L}^-).
\]
A simple argument using the Mayer-Vietoris sequence gives the
desired homology groups.
\eproof

Observe  that {\bf C4a} implies that $\dir(k)=-1$ for all $k>m$.
Therefore, one has the
following result.

\begin{cor}
 Assume
$
N\subset W\subset X_m
$
and the sequence
of pairs
$\{ a_k^\pm\in{\Bbb R}\mid a_k^-<a_k^+,\ k\in{\Bbb N}\}$ are
topologically self-consistent and satisfy {\bf C4a}. Fix an integer $r
>m$ and let
\[
\tilde{N}:= N\times [a_{m+1}^-,a_{m+1}^+]\times [a_{m+2}^-,a_{m+2}^+]
\times \cdots\times [a_{r}^-,a_{r}^+].
\]

Then
\[
CH_*(\Inv (\tilde{N})) \cong CH_*(\Inv (N)).
\]
\end{cor}

The following theorem  is used for all the results described in
the Introduction.

\begin{thm}
\label{thm:fixedpoint}
 Assume
$
N\subset W\subset X_m
$
and the sequence
of pairs
$\{ a_k^\pm\in{\Bbb R}\mid a_k^-<a_k^+,\ k\in{\Bbb N}\}$ are
topologically self-consistent and satisfy {\bf C4a}. Assume
\[
CH_j(\Inv(N,\varphi)) \cong\cases{ {\Bbb Z} & if $j=l$,\cr
0 & otherwise,}
\]
for some $l$, then there exists
\[
u^*\in N\times\prod_{k=m+1}^\infty [a_k^-,a_k^+],
\]
 a fixed point for the partial differential
equation (\ref{eq:pde}).
\end{thm}

\proof Combining Theorems~\ref{thm:indexcomp} and \ref{thm:fixedptindex},
immediately gives that for each $r>M$ there exists a fixed point
\begin{displaymath}
z_r\in N\times\prod_{k=m+1}^r [a_k^-,a_k^+]
\end{displaymath}
for the Galerkin projection onto the first $r$ coordinates.

Since $N\times\prod_{k=m+1}^\infty [a_k^-,a_k^+]$ is compact the
collection $\{ z_r\mid r = m+1,m+2,\ldots\}$ contains a limit
point $u^*$. From the continuity of $P_n\circ F$ on
$W \oplus \prod_{k=m+1}^\infty [a_k^-,a_k^+]$ it follows that
$P_n \circ F(u^*)=0$ for each $n \in {\Bbb N}$.
By Proposition~\ref{prop:ladder} $u^*$ is a fixed point for (\ref{eq:pde}).
\eproof

\subsection{Remarks on Related Work}

We are aware of at least two other results that are closely related
to the methods described in this Section. The first is work of
L. Cesari \cite{cesari} from the early 60's which in spirit is
very similar to ours. His method can be
characterized as follows \cite{williams}. Let $B$ be a Banach
space. Let $X$ be a finite dimensional subspace of $B$ and
let
$P:B\to X$ be a projection. Let $\tilde{N}\subset B$ be closed with
the property that $P\tilde{N}=N\subset X$ is compact and for every $x\in
N$,
$P^{-1}(x)\cap \tilde{N}$ is closed. Consider a continuous map
$f:\tilde{N}\to B$. The goal is to find fixed points for $f$ by studying
the behavior of the projection of the map onto $X$.

It is obvious that $u^*$ is a fixed point of $f$, if and only if
 $Pu^* = Pfu^*$, and
$u^* = Pu^*+(I-P)f(u^*)$.

Cesari's method applies if and only if the following three conditions
are satisfied:
\begin{enumerate}
\item[(i)] For each $x\in N$,
\[
P+(I-P)f:P^{-1}(x)\cap \tilde{N}\to P^{-1}(x)\cap \tilde{N}
\]
is a contraction.
\item[(ii)] Given condition (i), for each $x\in N$, there exists
a unique $u(x)\in \tilde{N}$ such that $Pu(x)+(I-P)f(u(x))=u(x)$. The
function $u:N \to \tilde{N}$ is continuous.
\item[(iii)] There are no fixed points of $Pfu:N\to N$ on the
boundary of $N$.
\end{enumerate}

Relating this back to the context of this paper, observe that
a fixed point for the partial differential equation is a fixed point
for any nonzero constant time  map of the corresponding
semi-flow. (iii) is closely related to the condition {\bf C5}. As stated
(ii) is not well defined unless (i) holds. {\bf C4} is the analogous
assumption to (i) and differs in two significant ways. A necessary
condition to have a contraction, is for the stronger assumption of {\bf
C4a} to hold. However, {\bf C4a} is not sufficient. An important point is
that we do not make any assumptions on the direction of the
vector field within $\tilde{N}$. Thus, condition {\bf C4a} is in
principle easier to verify than (i). On the other hand, this makes
it clear that we cannot guarantee uniqueness of the fixed point
given our assumptions.

%Observe that the Cesari method can be applied only to prove an existence of
%fixed points of (\ref{eq:pde}), while our method can be applied also to
%periodic solutions.

The other work is due to C.\ Conley and P. Fife \cite{conley-fife} and is
closely related to Theorem~\ref{thm:indexcomp}. Formulas of the
form (\ref{eq:indexshift}) are classical in the context of
product systems (see \cite{conley-cbms}). In \cite{conley-fife} one
finds a
similar formula, but in that context at the parameter value
for which one computes the index in the lower dimensional system,
there is no higher dimensional dynamics defined. However, the
higher dimensional system is defined for an arbitrarily small
perturbation. The key idea is that in the proper context the
lower dimensional system is normally hyperbolic. In this paper
we circumvent this type of assumption using isolating
blocks, {\bf C5}, and imposing {\bf C4}.

\section{Estimates for Kuramoto-Sivashinsky equation}

As the Hilbert space $H$ for the  Kuramoto-Sivashinsky equation (\ref{eq:KS})
we choose
the subspace of $L^2(-\pi,\pi)$ consisting of $2 \pi$-periodic and odd
functions.

Since $u(t,x)$ is odd its  Fourier expansion takes the form
\begin{equation}
\label{eq:FD}
u(t,x)=\sum_{k=-\infty}^{k=\infty} b_k(t) \exp(i kx)
\end{equation}
Since $u(t,x)$ is real, $b_k=\bar{b}_{-k}$. Substituting (\ref{eq:FD})
into (\ref{eq:KS}) gives the following equations
\begin{equation}
  {\dot b_k}=(k^2 - \nu k^4)b_k + \mbox{i}k \sum_{m=-\infty}^{m=\infty}
        b_m b_{k-m} \label{eq:KSf}
\end{equation}
Since we are interested in solutions with odd symmetry it follows
that $b_k$ are pure imaginary. Let
\[
a_k := \sqrt{-1}\,  b_k.
\]
Then, $a_k=-a_{-k}$ and $a_0=0$ which results in the following
infinite system of ordinary differential equations
\begin{equation}
\label{eq:KSfodd}
\dot{a}_k  =  k^2(1-\nu k^2)a_k -  k\sum_{n=1}^{k-1} a_n a_{k-n}
 +2k\sum_{n=1}^\infty
a_n a_{n+k}\quad k=1,2,3,\ldots
\end{equation}

We will use these equations to draw rigorous conclusions about
Kuramoto-Sivashinsky by finding self-consistent apriori bounds
for (\ref{eq:KS}) that satisfy the stronger condition {\bf C4a}
and then applying Theorem~\ref{thm:fixedpoint}.
To do this, however, we will need to understand the errors contributed
by ignoring the higher modes and the errors introduced by the use of
floating point arithmetic.

Let $m,M\in{\Bbb N}$ be fixed with $m\leq M$. Let $W\subset{\Bbb R}^m$
and  $\{ a_k^\pm\in{\Bbb R}\mid k\in{\Bbb N}\}$ satisfy conditions
{\bf C1} - {\bf C3} with the added constraints that
\[
W = \prod_{k=1}^m [a_k^-,a_k^+]
\]
and
\begin{equation}
\label{eq:powerdecay}
a_k^\pm = \pm \frac{C_s}{k^s},\quad k>M
\end{equation}
for some constant $C_s>0$ and integer $s >1$.

Though technically incorrect, it is perhaps useful for the reader
to think of the numerical approximation of the dynamics being computed
with respect to
the finite dimensional system
\begin{equation}
\label{eq:galerkin}
\dot{a}_k  =  k^2(1-\nu k^2)a_k -  k
\sum_{n=1}^{k-1} a_n a_{k-n}
 +2k\sum_{n=1}^{M-k}
a_n a_{n+k}\quad k=1,\ldots, m
\end{equation}
where $a_k = (a_k^-+a_k^+)/2$ for $k=m+1,\ldots M$. In doing so
it becomes clear that there are essentially three levels of
approximation that need to be dealt with. The first involves the terms
in the infinite tail $\{ a_k\mid k>M\}$. These are completely
absent from (\ref{eq:galerkin}) and therefore must be absorbed
as a fixed error term (think of the term $E(x)$ in
(\ref{eq:nonlinearode})). The power decay rule (\ref{eq:powerdecay})
will be used to determine this quantity. The
second, involves the terms
$\{ a_k\mid m<k\leq M\}$. In principle, one could set $m=M$, however
our strategy is to try to obtain better estimates for these terms
than can be expected by the general decay of (\ref{eq:powerdecay}).
However, these terms act as constants and hence can be viewed
as parameters for the system (\ref{eq:galerkin}). Finally, the
terms $\{ a_k\mid k=1,\ldots ,m\}$ are the actual variables for
the dynamical system being studied. It should also be kept in
mind that we need to lift the index information, and therefore
need to be able to verify {\bf C4a} for all $k$.

Of course, our goal is that of rigorous computations. Therefore each of
the above mentioned $a_i$ is actually an {\em interval}. The
intervals associated with $\{ a_k\mid k=1,\ldots ,m\}$ are
essentially determined by the floating point approximations. For
$k>m$, the intervals are
\[
a_k =[a_k^-,a_k^+].
\]

We will let
\begin{displaymath}
  |a_k|:=\max\{|a_k^-|, |a_k^+|\}.
\end{displaymath}

To compute the above mentioned errors we return to (\ref{eq:KSfodd})
and observe that in addition to the linear part there is a
finite sum of terms
\begin{equation}
\label{eq:FS}
FS(k) = \sum_{n=1}^{k-1} a_n a_{k-n}
\end{equation}
and an infinite sum of terms
\begin{equation}
\label{eq:IS}
IS(k) = \sum_{n=1}^\infty a_n a_{n+k}.
\end{equation}
Obviously bounds on these terms are necessary.

\subsection{$1\leq k\leq M$}

Since $FS(k)$ is a finite sum and we have already chosen the
interval values for $a_n$, we can explicitly compute $FS(k)$.
Perhaps it is worth noting that to evaluate $FS(k)$ only involves
the intervals $\{ a_n\mid n=1,\ldots, M-1\}$ which are chosen in
such a way that we expect them to be reasonably good approximations
of the actual terms.

\begin{lem}
\label{lem:ISm}
Assume $1\leq k\leq M$. Then,
\begin{eqnarray*}
IS(k) & \subset & \sum_{n=1}^{M-k}a_n a_{k+n} +
C_s\sum_{n=M-k+1}^M \frac{|a_n|}{(k+n)^s}[-1,1] + \\
& & \quad  \frac{C_s^2}{(k+M+1)^s(s-1)M^{s-1}}[-1,1]
\end{eqnarray*}
\end{lem}

\proof By definition,
\[
IS(k) = \sum_{n=1}^{M-k}a_n a_{k+n} + \sum_{n=M-k+1}^{M}a_n a_{k+n}
+ \sum_{n=M+1}^{\infty}a_n a_{k+n}.
\]
With regard to the second sum
\begin{eqnarray*}
 \sum_{n=M-k+1}^{M}a_n a_{k+n} &\subset &
 \sum_{n=M-k+1}^{M}|a_n| \frac{C_s}{(k+n)^s}[-1,1].
\end{eqnarray*}
Finally, the third sum produces
\begin{eqnarray*}
\sum_{n=M+1}^{\infty}a_n a_{k+n} &\subset &
\sum_{n=M+1}^{\infty}\frac{C_s}{n^s}\frac{C_s}{(n+k)^s}[-1,1] \\
 &\subset &
\frac{C_s^2}{(k+M+1)^s}[-1,1]\sum_{n=M+1}^{\infty}\frac{1}{n^s}
\\
 &\subset & \frac{C_s^2}{(k+M+1)^s(s-1)M^{s-1}}[-1,1] .
\end{eqnarray*}
In above derivation we used the following estimate
\begin{displaymath}
  \sum_{n=M+1}^\infty \frac{1}{n^s} < \int_{M}^\infty \frac{dx}{n^s}=
   \frac{1}{(s-1)M^{s-1}}
\end{displaymath}
\eproof

\begin{rem}
{\em
This estimate and some of those that follow  can be
improved by noting that
\[
\sum_{n=M+1}^{\infty}\frac{1}{n^s(n+k)^s} <
\int_M^\infty \frac{dx}{x^s(x+k)^s}.
\]
}
\end{rem}
Of course, the right hand side has an explicit rational expression,
but it is rather complicated for large $s$ and so was not utilized
here.

A simple extension of Lemma~\ref{lem:ISm} leads to the
following corollary.

\begin{cor}
\label{cor:tailest}
Let $1\leq k\leq m$. Then,
\begin{eqnarray*}
\sum_{n=m-k+1}^\infty a_n a_{n+k} &\subset &\sum_{n=m-k+1}^{M-k}a_n
a_{n+k} + C_s\sum_{n=M-k+1}^M \frac{|a_n|}{(k+n)^s}[-1,1] + \\
& & \quad  \frac{C_s^2}{(k+M+1)^s(s-1)M^{s-1}}[-1,1]
\end{eqnarray*}
\end{cor}

Observe that collorary \ref{cor:tailest}  estimates the error
in the vector field due to the Galerkin projection, namely
\begin{equation}
    A_k(p+q) - A_k F(p) =2k \sum_{n=m-k+1}^\infty a_n a_{n+k}
\end{equation}

\subsection{$k>M$}

Throughout this section it is assumed that $k>M$. Let
\[
e(k) :=\cases{1 & if $k$ is even,\cr
0 & if $k$ is odd.}
\]

\begin{lem}
\label{lem:M<k<2M}
Let $M<k\leq 2M$. Then,
\[
FS(k) \subset 2\sum_{n=k-M}^{\lfloor k/2\rfloor}a_n a_{k-n}+e(k)a_{k/2}^2
+2C_s\sum_{n=1}^{k-M-1}\frac{|a_n|}{(k-n)^s}[-1,1].
\]
\end{lem}
\proof Expanding (\ref{eq:FS}) gives
\begin{eqnarray*}
FS(k) & = &  \sum_{n=1}^{k-M-1} a_n a_{k-n} + \sum_{n=k-M}^M a_n
a_{k-n} +\sum_{n=M+1}^{k-1}a_n a_{k-n} \\
&\subset & 2C_s \sum_{n=1}^{k-M-1} \frac{|a_n|}{(k-n)^s}[-1,1] +
\sum_{n=k-M}^M a_n a_{k-n} \\
&\subset & \sum_{n=k-M}^{\lfloor k/2\rfloor}a_n a_{k-n}+e(k)a_{k/2}^2
+ 2C_s \sum_{n=1}^{k-M-1} \frac{|a_n|}{(k-n)^s}[-1,1].
\end{eqnarray*}
\eproof

\begin{lem}
\label{lem:k>2M}
Let $k>2M$. Then,
\[
FS(k) \subset \frac{C_s}{k^{s-1}}\left(\frac{2^{s+1}}{2M+1}\sum_{n=1}^M
|a_n| +\frac{C_s 4^s}{(2M+1)^{s+1}} + \frac{C_s 2^s}{(s-1)M^s}
\right)[-1,1].
\]
\end{lem}

\proof From (\ref{eq:FS}) it follows that
\begin{eqnarray*}
FS(k) & = & \sum_{n=1}^{k-1} a_n a_{k-n} \\
& = & 2\sum_{n=1}^{\lfloor k/2\rfloor}a_n a_{k-n}+e(k)a_{k/2}^2 \\
& = & 2\sum_{n=1}^{M}a_n a_{k-n} + 2\sum_{n=M+1}^{\lfloor k/2\rfloor}a_n
a_{k-n}+e(k)a_{k/2}^2.
\end{eqnarray*}
Each of these terms will be estimated separately. The first one results
in:
\begin{eqnarray*}
 \sum_{n=1}^{M}a_n a_{k-n} &\subset &
\sum_{n=1}^M \frac{|a_n|C_s}{(k-M)^s}[-1,1] \\
&\subset & \frac{C_s}{k^s(1-M/k)^s}[-1,1] \sum_{n=1}^M |a_n| \\
&\subset & \frac{2^{s}C_s}{k^s}[-1,1] \sum_{n=1}^M |a_n| \\
&\subset & \frac{2^{s}C_s}{k^{s-1}(2M+1)}[-1,1] \sum_{n=1}^M |a_n|.
\end{eqnarray*}
The second term leads to
\begin{eqnarray*}
\sum_{n=M+1}^{\lfloor k/2\rfloor}a_n a_{k-n} &\subset &
C_s^2\sum_{n=M+1}^{\lfloor k/2\rfloor}\frac{1}{n^s(k-n)^s}[-1,1]\\
&=& \frac{C_s^2}{k^s}\sum_{n=M+1}^{\lfloor
k/2\rfloor}\frac{1}{n^s(1-n/k)^s}[-1,1]\\
&\subset & \frac{C_s^2 2^s}{k^s}\sum_{n=M+1}^{\lfloor
k/2\rfloor}\frac{1}{n^s}[-1,1]\\
&\subset & \frac{C_s^2 2^s}{k^s}\int_{M}^\infty \frac{dx}{x}[-1,1]\\
&=& \frac{C_s^2 2^s}{k^s(s-1)M^{s-1}}[-1,1]\\
&\subset & \frac{C_s^2 2^{s-1}}{k^s(s-1)M^{s}}[-1,1].
\end{eqnarray*}
Finally, the third term gives rise to
\[
e(k)a_{k/2}^2\subset \frac{C_s^2 2^{2s}}{k^{2s}}[-1,1]
\subset \frac{1}{k^{s-1}}\frac{C_s^2 4^s}{(2M+1)^{s+1}}[-1,1].
\]
\eproof

Turning now to the infinite sum we can obtain the following estimate.

\begin{lem}
\label{lem:IS(k)est}
Let $k>M$. Then,
\[
IS(k) \subset \frac{C_s}{k^{s-1}(M+1)}\left(
\frac{C_s}{(M+1)^{s-1}(s-1)} + \sum_{n=1}^M |a_n|
\right)[-1,1].
\]
\end{lem}

\proof From (\ref{eq:IS}) it follows that
\[
IS(k) = \sum_{n=1}^M a_n a_{k+n} + \sum_{n=M+1}^\infty a_n a_{k+n}.
\]
As in the previous case, each term is treated separately.
\begin{eqnarray*}
\sum_{n=1}^M a_n a_{k+n} &\subset & C_s \sum_{n=1}^M
\frac{|a_n|}{(k+n)^s}[-1,1] \\
& = & \frac{C_s}{k^s}\sum_{n=1}^M \frac{|a_n|}{(1+n/k)^s}[-1,1] \\
 &\subset & \frac{C_s}{k^{s-1}(M+1)}\sum_{n=1}^M |a_n| [-1,1].
\end{eqnarray*}
The remaining term leads to
\begin{eqnarray*}
\sum_{n=M+1}^\infty a_n a_{k+n} & = & C_s^2 \sum_{n=M+1}^\infty
\frac{1}{m^s(k+m)^s}[-1,1] \\
&\subset & \frac{C_s^2}{(M+1)^s}\int_M^\infty \frac{dx}{(k+x)^s} [-1,1]\\
&=& \frac{C_s^2}{(M+1)^s(s-1)(k+M)^{s-1}} [-1,1]\\
&\subset & \frac{C_s^2}{(M+1)^s(s-1)k^{s-1}}[-1,1].
\end{eqnarray*}
\eproof

\subsection{Refining the Self-Consistent Bounds}
\label{sec:refinement}

The proof of our results obviously depends on having good
self-consistent  bounds and the precision of the final result
is determined directly by these bounds.  For this reason it
is important to have a process by which these bounds can be
improved. With this in mind consider an initial sequence of bounds
$\{ a_k^\pm\}$ which defines the sets
\[
W=\prod_{k=1}^m [a_k^- ,a_k^+]\quad{\rm and}\quad V=\prod_{k=m+1}^\infty
[a_k^- ,a_k^+].
\]
We will also assume that
\begin{displaymath}
  1 < \nu m^2.
\end{displaymath}
This condition means that the Fourier modes for $k \geq m$ are linearly stable.

We shall
describe the refinement procedure under the assumption of {\bf C4a}. In
particular, we need that our sequence $\{ a_k^\pm\}$ satisfy $\dir(k)=-1$
for all $k>m$. We also assume that since we can numerically solve
(\ref{eq:galode}), that the estimates for $W$ are reasonably good.

We will inductively adjust $a_k\pm$, for $k=m+1,\ldots ,M$, beginning with
$a_{m+1}^\pm$, as follows. Let $a\in W\oplus V$ such that
$a_k=a_k^+$. To satisfy {\bf C4a} requires that $\dot{a}_k<0$, i.e.
\[
k^2(1-\nu k^2)a_k^+ - kFS(k) + 2kIS(k) <0.
\]
This is equivalent to requiring
\begin{equation}
\label{eq:refininginequality}
a_k^+ > \frac{2IS(k) - FS(k)}{k^3(\nu - k^{-2})}.
\end{equation}
Of course, our goal is to make $a_k^+$ as small as possible within
the constraints imposed by the approximations. Since we are
iteratively improving our bounds, it is reasonable to assume that
a worst case equality is the best guess at this stage in the procedure.
Note, we are not claiming a proof at this point, we just are seeking
good bounds which later will be verified to be self-consistent bounds.
So using Lemma~\ref{lem:ISm}, define $f_k^{\pm}$ to be bounds
for $2IS(k) - FS(k)$,
\begin{eqnarray*}
  [f_k^-,f_k^+]:=2\sum_{n=1}^{M-k}a_n a_{k+n} +
2C_s\sum_{n=M-k+1}^M \frac{|a_n|}{(k+n)^s}[-1,1] + \\
  \frac{2C_s^2}{(k+M+1)^s(s-1)M^{s-1}}[-1,1]-
            \sum_{n=1}^{k-1} a_na_{k-n}.
\end{eqnarray*}

The new value of $a_k^+$ is given by
\[
a_k^+ := \frac{f_k^+}{k^3(\nu - k^{-2})}.
\]
A similar argument suggests setting
\[
a_k^- :=\frac{f_k^-}{k^3(\nu - k^{-2})}.
\]

This approach works up to $k=M$. Recall that for $k>M$, we set
$a_k^\pm =\pm C_s/k^s$. Here our goal is to improve the
power of convergence, i.e.\ we want to increase $s$. Again, since we are
trying to satisfy {\bf C4a} the basic
inequality which needs to be satisfied  is
(\ref{eq:refininginequality}). The estimates for $FS(k)$ and
$IS(k)$ for $k>M$ obviously are crucial here. However, we had
two sets of estimates one for $M<k\leq 2M $ and the other for
$k>2M$. Thus, we need to choose the worst of
both estimates. This is done as follows.

Given an interval
$I\subset {\Bbb R}$ let
\[
|I| := \sup_{x\in I}|x|.
\]
With the estimate from Lemma~\ref{lem:M<k<2M}  in mind define
\[
D_1(k) :=\left| 2\sum_{n=k-M}^{\lfloor k/2\rfloor}a_n
a_{k-n}+e(k)a_{k/2}^2\right| +2C_s\sum_{n=1}^{k-M-1}\frac{|a_n|}{(k-n)^s}.
\]
Combining this with the estimate on $IS(k)$ given by
Lemma~\ref{lem:IS(k)est} and multiplying by $k^{s-1}$ leads to the
following definition
\[
D_1 := \frac{2C_s}{(M+1)}\left(
\frac{2C_s}{(M+1)^{s-1}(s-1)} + \sum_{n=1}^M |a_n|
\right) +\max_{M<k\leq 2M}k^{s-1}D_1(k).
\]

Turning now to the bounds for $k>2M$, Lemmas~\ref{lem:k>2M},
\ref{lem:IS(k)est}, and again multiplying by $k^{s-1}$ suggests
setting
\begin{eqnarray*}
D_2 &:=& \frac{2C_s}{(M+1)}\left(
\frac{2C_s}{(M+1)^{s-1}(s-1)} + \sum_{n=1}^M |a_n|
\right) + \\
&&\quad\quad
\frac{C_s}{k^{s-1}}\left(\frac{2^{s+1}}{2M+1}\sum_{n=1}^M
|a_n| +\frac{C_s 4^s}{(2M+1)^{s+1}} + \frac{C_s 2^s}{(s-1)M^s}.
\right)
\end{eqnarray*}

From Lemmas~\ref{lem:M<k<2M}, \ref{lem:k>2M}, and \ref{lem:IS(k)est}
we obtain the following result.

\begin{cor}
For $k>M$ and $D_s:=\max \{ D_1,D_2\}$
\[
\left| -FS(k)+2IS(k)\right | < \frac{D_s}{k^{s-1}}.
\]
\end{cor}

We will use this corollary to update the decay rate for the tail terms.
Again, we want (\ref{eq:refininginequality}) (which gave us {\bf C4a})
 to hold for all $k>M$. It is sufficient that
\[
a_k^+ >\frac{\left| -FS(k)+2IS(k)\right |}{k^3(\nu-(M+1)^{-2})}
\]
and therefore it is sufficient that
\[
a_k^+ >\frac{D_s}{k^{s-1}}\cdot\frac{1}{k^3}\cdot\frac{1}{\nu-(M+1)^{-2}}
=\frac{1}{k^{s+2}}\cdot\frac{D_s}{\nu-(M+1)^{-2}}.
\]
There is a similar inequality for $a_k^-$.

Setting this to an equality we can define
\begin{equation}
a_k^\pm := \pm\frac{C_{s+2}}{k^{s+2}}\quad{\rm where}\quad
C_{s+2}:=\frac{D_s}{\nu-(M+1)^{-2}}
\end{equation}
for $k>M$.

\section{A Typical Proof}

This section  describes the proof of the following
result.

\begin{thm}
\label{thm:proof}
 Let
\[
u(x) = \frac{1}{\sqrt{2}}\sin x - \frac{1}{8} \sin 2x.
\]
For $\nu=0.75$ there
exists an equilibrium solution $u^*(x)$ to (\ref{eq:KS})
such that
\[
|| u^* - u||_{L^2} < 0.052\quad{\rm  and}\quad
||u^*-u||_{C^0} < 0.05
\]
\end{thm}

\begin{figure}[h]
\begin{center}
\includegraphics[width=0.80\textwidth]{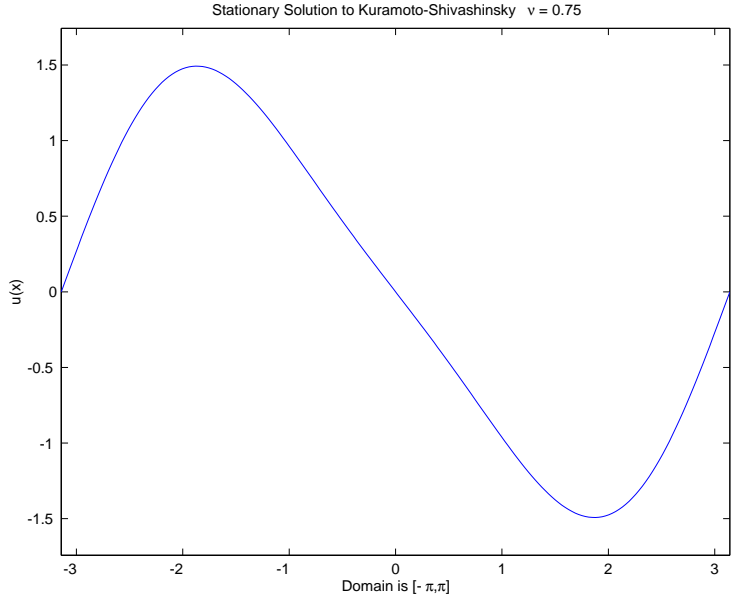}
\end{center}
\caption{The computed function $u(x)$ from Theorem \ref{thm:proof}}
\label{fig:equil2}
\end{figure}

The reader should observe that this is a weaker version of
Theorem~\ref{thm:stable}. However, we present its proof since it contains
all the essential features, but with a very low dimensional approximation.
The rest of the results described in the introduction were proved in
a similar manner.

The first step is to choose $m$ the dimension of the Galerkin
approximation and $M$ the level to which we make specific choices
for the $\{ a_k^\pm\}$. For $m=2$, (\ref{eq:galerkin}) reduces to the
system of equations
\begin{eqnarray}
\label{eq:ksgal22}
\dot{a}_1& = &\frac{1}{4}a_1 + 2a_1a_2   \\
\dot{a}_2&= &-8a_2-2a_1^2.  \nonumber
\end{eqnarray}
A simple algebraic computation shows that this system has
exactly three fixed points: $(0,0)$-unstable, with one-dimensional
unstable manifold and two attracting fixed points $u_{\pm}=(\pm
\frac{1}{\sqrt{2}},-\frac{1}{8})$. Theorem~\ref{thm:proof} is
obtained by studying the dynamics of (\ref{eq:ksgal22}) in a neighborhood
of $(\frac{1}{\sqrt{2}},-\frac{1}{8})$.

The next step is to obtain self-consistent apriori bounds for
(\ref{eq:KSfodd}). It is unrealistic to expect that goods bounds
can be obtained immediately. Thus, we make a reasonable guess
for bounds and then try to improve them. Let
\[
W= [\frac{1}{\sqrt{2}}-0.1,\frac{1}{\sqrt{2}}+0.1]
\times [-\frac{1}{8}-0.1,-\frac{1}{8}+0.1].
\]
The initial estimates for $[a_k^-,a_k^+]$ are given in
Table~\ref{table:init}.  The formula used to derive these initial
estimates will be presented in Section~\ref{sec:initialest}. The reason
for delaying the presentation is to emphasize the fact that the initial
estimates are only estimates. Obviously, choosing good estimates
allows for faster convergence and choosing terrible estimates
will probably result in a failure of convergence.

\begin{table}[h]
\begin{center}
\begin{tabular}{|r|c|}
      \hline
       k & initial $[a^-_k,a^+_k]$
\\
       \hline
3 & [-0.157542,\ 0.157542] \\
4 & [-0.0463226,\ 0.0463226]  \\
5 & [-0.0183725,\ 0.0183725]   \\
6 & [-0.00871023,\ 0.00871023]   \\
7 & [-0.00465407,\ 0.00465407]  \\
8 & [-0.00271036,\ 0.00271036]   \\
9 & [-0.00168454,\ 0.00168454]    \\
10 & [-0.00110173,\ 0.00110173]    \\
$>10$ & $[-1,1]\cdot 10.9915/k^4$   \\
   \hline
    \end{tabular}
\end{center}
\caption{Initial estimates for the intervals $[a_k^-,a_k^+]$.}
\label{table:init}
\end{table}

Beginning with the data in Table~\ref{table:init}, the refinement
procedure described in Section~\ref{sec:refinement} is used to update
$a_k^\pm$ for $k>2$. After three iterations one  obtains the  estimates
given in Table~\ref{table:final}. It can now be checked that  $W$
and $\{ a_k^\pm\}$ from Table~\ref{table:final} form self-consistent
bounds.

\begin{table}[h]
\begin{center}
\begin{tabular}{|r|c|}
\hline
k & final $[a^-_k,a^+_k]$ \\
       \hline
3 & $ [0,\ 0.021055] $ \\
4 & $ [-0.00192301,\ 0] $\\
5 & $ [-1.8253\times 10^{-07},\ 0.000141734] $  \\
6 & $ [-9.85549\times 10^{-06},\ 8.64999\times 10^{-09}] $ \\
7 & $ [-6.55526\times 10^{-10},\ 6.42034\times 10^{-07}] $\\
8 & $ [-4.03088\times 10^{-08},\ 9.30992\times 10^{-11}] $\\
9 & $ [-3.51558\times 10^{-10},\ 2.79203\times 10^{-09}] $  \\
10 & $ [-1.11597\times 10^{-09},\ 9.71368\times 10^{-10}] $ \\
$>10$ & $[-1,1]\cdot 10285.3/k^{10}$  \\
   \hline
    \end{tabular}
\end{center}
\caption{Estimates for the intervals $[a_k^-,a_k^+]$  representing
 self-consistent apriori bounds}
\label{table:final}
\end{table}

What should be clear at this point is that the uncertainty contributed
by the terms not in the Galerkin projection are extremely small.
Obviously, at this point most of the uncertainty
is due to the size of $W$.

Having controlled the errors from the Galerkin truncation, the next
step is to obtain an isolating block $N$ which is topologically
self consistent with $W$ and $\{ a_k\}$. In determing $N$ we
use the vector field (\ref{eq:ksgal22}). Of course, the correct
equations are given by (\ref{eq:KSfodd}):
\begin{eqnarray*}
\dot{a}_1 & = & \frac{1}{4}a_1 +2\sum_{n=1}^\infty a_n a_{n+1}\\
\dot{a}_2 & = & -8a_k -  2a_1^2 + 4\sum_{n=1}^\infty a_n a_{n+2}.
\end{eqnarray*}
Using Corollary~\ref{cor:tailest} one obtains the following bounds on the
errors, $\epsilon_i$, $i=1,2$,
\[
\epsilon_1=[-0.00955626,7.06005 \cdot 10^{-10}] \quad{\rm and}\quad
\epsilon_2=[-1.544 \cdot 10^{-8},0.0697171].
\]

Thus, the equations for which the isolating neighborhood should be
found are
\begin{eqnarray}
\label{eq:ksgal22err}
\dot{a}_1  & = &\frac{1}{4}a_1 + 2a_1a_2 +\epsilon_1 \\
\dot{a}_2 &= &-8a_2-2a_1^2 +\epsilon_2.  \nonumber
\end{eqnarray}
The construction of the isolating block around $(1/\sqrt{2},-1/8)$
is easier if one works in coordinates determined by the
eigenfunctions of the linearized equations at the
fixed point. The eigenvalues are
\[
\lambda_1=-4 + 2\sqrt{3}, \quad \lambda_2=-4 - 2\sqrt{3}
\]
with corresponding unit eigenvectors
\[
v_1=\frac{1}{\sqrt{15-8\sqrt{3}}}(1,\frac{\lambda_1}{\sqrt{2}})^t, \quad
v_2=\frac{1}{\sqrt{15+8\sqrt{3}}}(1,\frac{\lambda_2}{\sqrt{2}})^t.
\]
Let $T$ be the affine change of variables from $(x,y)^t$ in this
new basis to the original variables $(a_1,a_2)^t$. Then, on the
set $T^{-1}(W)$ (\ref{eq:ksgal22err}) becomes
\begin{eqnarray}
\label{eq:ksg2t1}
\dot{x}=(-4 + 2\sqrt{3}) x + f_1(x,y) + {\tilde \epsilon_1}   \\
\dot{y}=(-4 - 2\sqrt{3}) y + f_2(x,y) + {\tilde \epsilon_2} \nonumber
\end{eqnarray}
where $f_1$ and $f_2$ are polynomials containing only terms of degree two
and ${\tilde \epsilon_1}$, ${\tilde \epsilon_2}$ are obtained from
$\epsilon_1$,
$\epsilon_2$ by the transformation $T$. In particular,
\[
{\tilde \epsilon_1}=[-0.0110098,0.0152184]  ,
\quad {\tilde \epsilon_2}=[-0.0764461,0.00397075]
\]

Set ${\tilde W}=[-0.0748016,0.0748016]^2$, then $T({\tilde W}) \subset W$.
Thus, an isolating block $\tilde{N}\subset\tilde{W}$ (which satisfies
the error constraints) will
give rise to an isolating block $N=T(\tilde{N})\subset W$ such that
$N$, $W$ and $\{ a_k\}$ are topologically consistent. Since for each
$k$, $\dir(k)=-1$,
\[
CH_j(\Inv (N))\cong\cases{ {\Bbb Z} & if $j=0$,\cr
0 & otherwise.}
\]
Thus, by Theorem~\ref{thm:fixedpoint}, there exists the desired fixed
point in
\[
u^*\in N\times \prod_{k=2}^\infty [a_k^-,a_k^+].
\]

Thus, all that remains is to construct $\tilde{N}$. In the $(x,y)$
coordinates the error constraints become
\begin{eqnarray*}
  f_1(x,y) + {\tilde \epsilon_1} & \subset &
(b_x,B_x):=(-0.0334071,0.0376157)
   \quad \mbox{for $(x,y) \in {\tilde W}$} \\
  f_2(x,y) + {\tilde \epsilon_2} & \subset &
(b_y,B_y):=(-0.0764461,0.00397075)
   \quad \mbox{for $(x,y) \in {\tilde W}$} \\
\end{eqnarray*}
This implies that bounds on the derivative are given by
\[
 \lambda_1 (x + \frac{b_x}{\lambda_1})  < \dot{x}  <
             \lambda_1 (x + \frac{B_x}{\lambda_1}) \quad{\rm and}\quad
    \lambda_2 (y + \frac{b_y}{\lambda_2})<  \dot{y}  <
                    \lambda_2 (y + \frac{B_y}{\lambda_2}).
\]

Observe that $\lambda_i < 0$, hence the box
\begin{eqnarray*}
  {\tilde N} & = &
\left[-\frac{b_x}{\lambda_1},-\frac{B_x}{\lambda_1}\right]
\times
   \left[-\frac{b_y}{\lambda_2},-\frac{B_y}{\lambda_2}\right]  \\
& = & {[}-0.0623385, 0.0701918{]} \times [-0.0132425,0.00353264] \subset
{\tilde W}
\end{eqnarray*}
is an isolating block.

\section{Obtaining the Initial Estimates}
\label{sec:initialest}

Before beginning this section we want to once again emphasize
that the proofs of the theorems in this paper are in principle
independent of this section. On the other hand, good initial
guesses greatly improve the speed of convergence. The estimates
described in what follows apparently provide excellent initial
values for the self-consistent bounds.

We will follow  the arguments from \cite{jkt} to produce  estimates
for errors in the Galerkin projection. 
Kuramoto-Sivashinsky can be written in the form
\[
u_t + \nu Au - A^{1/2}u + 2 B(u,u)=0  
\]
where
\[
A= \frac{\partial^4}{\partial x^4}, \quad
A^{1/2}=-\frac{\partial^2}{\partial x^2}, \quad
B(u,v)=u \frac{\partial v }{\partial x}\ .
\]
While $A^{1/4} \neq \frac{\partial}{\partial x}$, it is still the case
that
\[
|A^{1/4}u|_2=|\frac{\partial u}{\partial x}|_2 .
\]
To simplify the notation, let
\[
 |u|=|u|_2, \quad  \|u\|=|A^{1/4}u|_2.
\]

Since we are interested in bounded invariant sets we can without loss
of generality assume  the following apriori bounds for the invariant set
under consideration:
\begin{equation}
  |u(t)| \leq \rho_0, \quad \|u(t)\| \leq \rho_1 , \quad \mbox{for $t > T(u(0))$}
  \label{eq:abounds}
\end{equation}

We will make use of the following inequality \cite[Lemma 1.4]{jkt}
\begin{equation}
\label{eq:1.4}
|(B(u,v),w)| \leq \sqrt{2}|u|^{1/2}\|u\|^{1/2}\|v\||w|
\end{equation}

The eigenvalues of $A$ are $\lambda_n=n^4$, $n \in {\Bbb N}$ and 
the corresponding complete family of orthonormal 
eigenfunctions are $\{ \frac{1}{\sqrt{\pi}}\sin(nx)\}$.
Let $P=P_m$ be an orthogonal projection on first  $m$  eigenfunctions and
set $Q=I - P$.

Using the  decomposition
\[
  p=Pu, \quad q = Qu
\]
and the same abuse of notation
the equation for $q$ is
\begin{equation}
\label{eq:qdot}
\dot{q}= - \nu Aq +A^{1/2}q - 2 Q B(u,u).
\end{equation}
\begin{thm}
\label{th:gale1}
Under the assumptions stated above
if $m$ is large enough such that $\lambda_{m+1} > \frac{1}{\nu^2}$, then
\[
  \limsup_{t \to \infty}|q(t)| \leq \frac{2 \sqrt{2} \rho_0^{1/2}
\rho_1^{3/2}}
    {\lambda_{m+1}(\nu - \lambda_{m+1}^{-\frac{1}{2}})}.
\]
\end{thm}
\noindent
{\bf Proof:} Beginning with (\ref{eq:qdot}) and  taking a scalar product
with $q$ gives
\[
  (\frac{dq}{dt}|q) = -\nu (Aq|q)+(A^{1/2}q|q) - 2 (B(u,u)|q).
\]
Therefore,
\begin{equation}
 \frac{1}{2} \frac{d}{dt}|q|^2 \leq  -\nu (Aq|q) +(A^{1/2}q|q) + 2
|(B(u,u)|q)|
  \label{eq:aimgi}
\end{equation}
Observe that
\begin{equation}
  (Aq|q)=(A^{1/2}q|A^{1/2}q)=|A^{1/2}q|^2 \label{eq:aqq}
\end{equation}
and
\begin{eqnarray}
(A^{1/2}q|q) & = &\sum_{n=m+1}^{\infty} \lambda_n^{1/2}|q_n|^2
\nonumber\\ 
& = & \lambda_{m+1}^{-1/2} \sum_{n=m+1}^{\infty} \lambda_{m+1}^{1/2}
               \lambda_n^{1/2}|q_n|^2    \nonumber \\
&\leq &   \lambda_{m+1}^{-1/2} \sum_{n=m+1}^{\infty} \lambda_n|q_n|^2 
\nonumber \\
& = &  \lambda_{m+1}^{-1/2}|A^{1/2}q|^2 \label{eq:asqq}
\end{eqnarray}
Thus, (\ref{eq:aimgi}) (\ref{eq:aqq}) and (\ref{eq:asqq}) imply that
\begin{equation}
\label{eq:inter}
 \frac{d}{dt}|q|^2 \leq - 2\nu |A^{1/2}q|^2 +
2\lambda_{m+1}^{-1/2}|A^{1/2}q|^2
   + 4 |(B(u,u)|q)|.
\end{equation}

From (\ref{eq:1.4}) it follows that
\[
    |(B(u,u),q)| \leq \sqrt{2}|u|^{1/2}\|u\|^{3/2}|q|.
\]
While, (\ref{eq:abounds}) implies that for $t > T(u)$
\[
  |(B(u,u),q)| \leq \sqrt{2}\rho_0^{1/2}\rho_1^{3/2}|q|
\]
Therefore, (\ref{eq:inter}) becomes
\[
 \frac{d}{dt}|q|^2 \leq  -2(\nu - \lambda_{m+1}^{-1/2})|A^{1/2}q|^2 +
       4 \sqrt{2}\rho_0^{1/2}\rho_1^{3/2}|q|
\]

Observe that  $|A^{1/2}q|^2 \geq \lambda_{m+1} |q|^2$ and by 
assumption
$\nu - \lambda_{m+1}^{-1/2} >0$. Thus,
\begin{eqnarray*}
 \frac{d}{dt}|q|^2 &\leq & - 2(\nu -
\lambda_{m+1}^{-1/2})\lambda_{m+1}|q|^2 
        +4\sqrt{2}\rho_0^{1/2}\rho_1^{3/2}|q| \\
& = &(4\sqrt{2}\rho_0^{1/2}\rho_1^{3/2}- 2(\nu
-\lambda_{m+1}^{-1/2})\lambda_{m+1}|q|)|q|
\end{eqnarray*}
Thus, for 
\[
|q| > \frac{4\sqrt{2}\rho_0^{1/2}\rho_1^{3/2}}{2(\nu
-\lambda_{m+1}^{-1/2})\lambda_{m+1}}
\]
 $\frac{d}{dt}|q|^2 < 0$ and hence,
\[
  \limsup_{t->\infty}|q(t)| \leq
\frac{4\sqrt{2}\rho_0^{1/2}\rho_1^{3/2}}{2(\nu
-\lambda_{m+1}^{-1/2})\lambda_{m+1}}.
\]
\eproof

\begin{rem}
{\em
Because an orthonormal collection eigenvectors were used for
the calculations in this section, the coefficients $q_k$ and
$a_k$ differ by a scaling, i.e.
\[
q_k = -2\sqrt{\pi}a_n.
\]
}
\end{rem}

Theorem~\ref{th:gale1} can be used as follows. For a fixed $\nu$,
numerical experiments can suggest values for $\rho_0$ and $\rho_1$.
For $m<k\leq M$ one can use the formula
\[
a_k^\pm := \min\left\{ \pm\rho_0,\pm\frac{\rho_1}{k},
\pm\frac{4\sqrt{2\pi\rho_0\rho_1^3}}{k^4(\nu -k^-2)}\right\}.
\]
For $k>M$, one defines
\[
C_4:= \frac{4\sqrt{2\pi\rho_0\rho_1^3}}{\nu -(M+1)^{-2}}.
\]


\begin{thebibliography}{1}

\bibitem{az} G.\ Arioli and P. Zgliczy\'nski,
    Symbolic dynamics for the H\'enon--Heiles hamiltonian
       on the critical energy level, {\em J. Diff. Eq}, accepted

\bibitem{cime} L.\ Arnold, et. al., {\em Dynamical Systems},
Lect.\ Notes Math. 1609, R. Johnson, ed., 1995.

\bibitem{cesari} L.\ Cesari, Functional analysis and Galerkin's method,
{\em Mich.\ Math.\ Jour.} {\bf 11} (1964) 383-414.

\bibitem{conley-cbms} C.\ Conley,  {\em Isolated Invariant Sets and the Morse
Index.} CBMS Lecture Notes {\bf 38} A.M.S. Providence, R.I. 1978.

\bibitem{conley-fife} C.\ Conley and P. Fife, Critical manifolds,
travelling waves and an example from population genetics, {\em
J.\ Math.\ Bio.} {\bf 14} (1982) 159-176.

\bibitem{FNST} C. Foias, B. Nicolaenko, G. Sell, R. Temam,
   Inertial manifolds for the Kuramoto--Sivashinsky equation
    and an estimate of their lowest dimension.
    {\em J. Math. Pures Appl.} {\bf 67}, (1988), 197--226

\bibitem{gz} Z.\ Galias and P. Zgliczy\'nski,
Computer assisted proof of chaos in the Lorenz
  system,  {\em Physica D}, {\bf 115}, 1998,165--188


\bibitem{hale} J.~K. Hale, {\em Asymptotic Behavior of Dissipative
Systems}, Math. Surveys and Monographs {\bf 25}, AMS 1988


\bibitem{hmo} J.~K. Hale, L. T. Magalh\~aes,  and W. M. Oliva,
{\em An Introduction to Infinite Dimensional Dynamical Systems -
Geometric Theory}, Appl.\ Math.\ Sci.\ {\bf 47}, Springer-Verlag 1984.

\bibitem{jkt} M. Jolly, I. Kevrekidis, E. Titi,  Approximate
    inertial manifolds for the Kuramoto--Sivashinsky equation: analysis
    and computations, {\em Physica D} {\bf 44}, 38-60, (1990)


\bibitem{JRT} M. Jolly, R. Rosa, R. Temam,  Evaluating the dimension
   of an inertial manifold for the Kuramoto--Sivashinsky  Equation, preprint

\bibitem{KT} Y. Kuramoto, T. Tsuzuki,
Persistent propagation of concentration waves in dissipative media
   far from thermal equilibrium, {\em Prog. Theor. Phys.} {\bf 55},(1976), 365



\bibitem{mccord} C. McCord, Mappings and homological properties in the
Conley index theory, {\em Ergod.\ Th.\ \& Dyn.\ Sys.} {\bf 8*}
(1988) 175-198.

\bibitem{mm-lorenz} K.\ Mischaikow and M.\ Mrozek, Isolating neighborhoods
and Chaos,
     {\it Jap. J. Ind. \& Appl. Math.}, {\bf 12}, 1995, 205-236..

\bibitem{mm-lorenz2} K.\ Mischaikow and M.\ Mrozek,
      Chaos in Lorenz equations: a computer assisted proof,
      {\it Bull. Amer. Math. Soc. (N.S.)}, {\bf 33}(1995), 66-72.

\bibitem{handbook} K. Mischaikow and M. Mrozek, Conley Index Theory,
to appear in Handbook of Dynamical Systems, Vol. 3.

\bibitem{mms} K.\ Mischaikow, M.\ Mrozek and A.\ Szymczak, Chaos in
Lorenz  equations: a computer assisted proof, Part III:  Classical parameter
  values,
{\em JDE} to appear.

\bibitem{salamon} D.\ Salamon,  Connected simple systems and the Conley index
of isolated invariant sets.  {\em Trans. A. M. S.} {\bf 291} (1985) 1 - 41.

\bibitem{Si} G.I. Sivashinsky,
  Nonlinear analysis of hydrodynamical instability in laminar flames --
   1. Derivation of basic equations,
   {\em Acta Astron.} {\bf 4}, (1977), 1177


\bibitem{smoller} J.\ Smoller, {\em Shock Waves and Reaction Diffusion Equations},
Springer Verlag, New York, 1980.

\bibitem{temam} R. Temam, {\em Infinite-Dimensional Dynamical
Systems in Mechanics and Physics}, Springer-Verlag 1988.

\bibitem{williams} S. Williams, {\em On the Cesari Fixed Point
Method in a Banach Space}, Ph.D. Thesis, California Institute of
Technology, 1967.

\bibitem{wz} K.\ W\'ojcik and P.\ Zgliczy\'nski,
  Isolating segments, fixed point index and symbolic dynamics,
    {\em J. Diff. Eq}, 161, 245--288, (2000)

\bibitem{zg1} P.\ Zgliczy\'nski, Computer assisted proof of chaos
   in  the  H\'enon map and in the R\"ossler equations,
         {\em Nonlinearity},{\bf 10}, 1997,  No. 1, 243--252

\bibitem{zg2}  P.\ Zgliczy\'nski,
    Multidimensional perturbations of one-dimensional maps
           and  stability of Sharkovskii ordering,
       {\em Int. J. of Bifurcation and Chaos}, {\bf 9}, No. 9, 1999, 1867--1876

\end{thebibliography}
\end{document}